\documentclass[12pt, 
]{amsart}

\usepackage{amssymb}
\usepackage{amsmath}
\usepackage{amsthm}
\usepackage{amscd}
\usepackage{graphicx}
\usepackage[all]{xy}

\numberwithin{equation}{section}
\newtheorem{theo}{Theorem}[section]
\newtheorem{prop}[theo]{Proposition}

\newtheorem{cor}[theo]{Corollary}

\newtheorem{prob}[theo]{Problem}
\newtheorem{conj}[theo]{Conjecture}

\theoremstyle{definition}

\theoremstyle{remark}
\newtheorem{rem}[theo]{Remark}
\theoremstyle{definition}

\newcommand{\Supp}[0]{{\operatorname{Supp}}}

\newcommand{\ddbar}{dd^c}

\newcommand{\dbar}{\overline{\partial}}
\newcommand{\e}{\varepsilon}
\newcommand{\ome}{\widetilde{\omega}}
\newcommand{\I}[1]{\mathcal{I}(#1)}

\newcommand{\lla}[0]{{\langle\!\langle}}
\newcommand{\rra}[0]{{\rangle\!\rangle}}
\begin{document}

\title[]{Injectivity theorems with multiplier ideal sheaves \\and their applications}

\author{SHIN-ICHI MATSUMURA}

\address{Mathematical Institute, Tohoku University, 6-3, Aramaki Aza-Aoba, Aoba-ku, Sendai 980-8578, Japan. }

\email{{\tt
mshinichi@m.tohoku.ac.jp, mshinichi0@gmail.com}}

\thanks{Classification AMS 2010: 14F18, 32L10, 32L20. }


\maketitle

\begin{abstract}
The purpose of this survey is to present 
analytic versions of the injectivity theorem 
and their applications. 
The proof of our injectivity theorems is based on a 
combination of the $L^2$-method for the $\dbar$-equation and 
the theory of harmonic integrals.   
As applications, 
we obtain Nadel type vanishing theorems and 
extension theorems for pluri-canonical sections of log pairs. 
Moreover, we give some results on semi-ampleness 
related to the abundance conjecture in birational geometry 
(the minimal model program). 
\end{abstract}

\section{Introduction}\label{Intro}

The Kodaira vanishing theorem is 
one of the most celebrated results in complex geometry, 
and such results play an important role 
when we consider certain fundamental
problems in algebraic geometry and 
the theory of several complex variables, 
including asymptotics of linear systems, 
extension problems of holomorphic sections, the
minimal model program, and so on.  
According to these objectives,
the study of vanishing theorems has been 
a constant subject of attention in the last decades.
This paper is a survey of recent results in 
\cite{Mat13b} and \cite{GM13},  
whose purpose is to present generalizations of 
the Kodaira vanishing to pseudo-effective line bundles 
equipped with singular metrics and their applications, 
from the viewpoint of the theory of several complex variables and differential geometry.

\subsection{Analytic versions of the injectivity theorem}
In this subsection, we introduce analytic versions of 
the injectivity theorem. 
The injectivity theorem is 
a generalization of the vanishing theorem to 
\lq \lq semi-positive" line bundles, 
and it has been studied by several authors, for example, Tankeev (\cite{Tan71}), Koll\'ar (\cite{Kol86}), 
Enoki (\cite{Eno90}), Esnault-Viehweg (\cite{EV92}), 
Ohsawa (\cite{Ohs04}), Fujino (\cite{Fuj12a}, \cite{Fuj12b}), 
Ambro (\cite{Amb03}, \cite{Amb12}), and so on. 
In his paper \cite{Kol86}, Koll$\rm{\acute{a}}$r 
gave the following injectivity theorem for semi-ample line bundles, 
whose proof depends on the Hodge theory. 
In \cite{Eno90}, 
Enoki relaxed his assumption 
by a different method depending on the theory of harmonic integrals, 
which enables us to approach the injectivity theorem 
from the viewpoint of complex differential geometry.

\begin{theo}
[{\cite{Kol86} (resp. \cite{Eno90})}] 
\label{Kol}
Let $F$ be a semi-ample $($resp. semi-positive$)$ 
line bundle on a smooth projective variety 
$($resp. a compact K\"ahler manifold$)$ $X$. 
Then for a $($non-zero$)$ section $s$ of a positive multiple 
$F^{m}$ of the line bundle $F$, 
the multiplication map 
induced by the tensor product with $s$ 
\begin{equation*}
\Phi_{s}: H^{q}(X, K_{X} \otimes F) 
\xrightarrow{\otimes s} 
H^{q}(X, K_{X} \otimes F^{m+1} )
\end{equation*}
is injective for any $q$. 
Here $K_{X}$ denotes the canonical bundle of $X$. 
\end{theo}

The above theorem can be regarded as 
a generalization of the Kodaira vanishing theorem 
to semi-ample (semi-positive) line bundles. 
On the other hand, 
the Kodaira vanishing theorem has 
been generalized by Nadel (\cite{Nad89}, \cite{Nad90}). 
This generalization uses singular metrics with 
positive curvature and corresponds to 
the Kawamata-Viehweg vanishing theorem 
in algebraic geometry (\cite{Kaw82}, \cite{Vie82}). 
Therefore, in the same direction as this generalization, 
it is natural and of interest 
to study injectivity theorems for line bundles 
equipped with \lq \lq singular metrics".

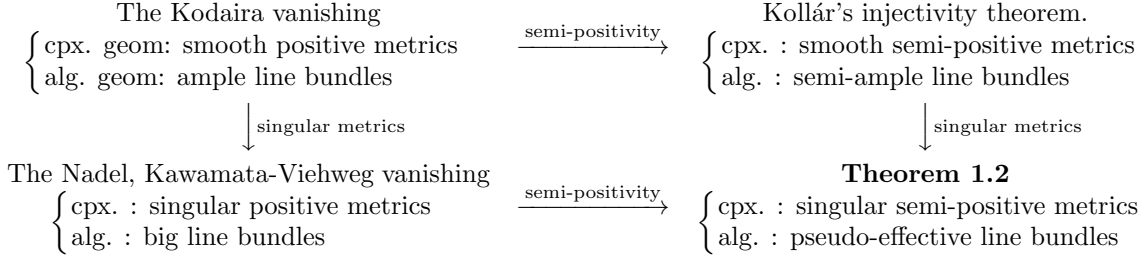
$$\hspace{-1.0cm} 
\footnotesize{
\begin{CD}
\begin{array}{cl}
\text{The Kodaira vanishing} \\
\left \{\hspace{-0.5cm} 
\begin{array}{cl}
&\text{cpx. geom: smooth positive metrics}\\
&\text{alg. geom: ample line bundles} 
\end{array}
\right.
\end{array}
 @>\text{semi-positivity}>> 
\begin{array}{cl}
\text{Koll\'ar's injectivity theorem.}\\
\left \{\hspace{-0.5cm} 
\begin{array}{cl}
&\text{cpx. : smooth semi-positive metrics}  \\
&\text{alg. : semi-ample line bundles} 
\end{array}
\right.
\end{array}
\\
@VV\text{singular metrics}V @VV\text{singular metrics}V \\
\begin{array}{cl}
\text{The Nadel, Kawamata-Viehweg vanishing}\\
\left \{\hspace{-0.5cm} 
\begin{array}{cl}
&\text{cpx. : singular positive metrics}\\
&\text{alg. : big line bundles} 
\end{array}
\right.
\end{array} 
@>\text{semi-positivity}>> 
\begin{array}{cl}
\text{{\bf{Theorem \ref{main-inj}}}}\\
\left \{\hspace{-0.5cm}
\begin{array}{cl}
&\text{cpx. : singular semi-positive metrics}\\
&\text{alg. : pseudo-effective line bundles} 
\end{array}
\right. 
\end{array}  \\
\end{CD}
}
$$
\vspace{0.1cm}
\\
The following theorem is one of the main results, 
which can be seen as a generalization of the injectivity theorem and the Nadel vanishing theorem.

\begin{theo}[\text{\cite[Theorem 1.3]{Mat13b}}]\label{main-inj}
Let $F$ be a line bundle on a compact K\"ahler manifold $X$ and 
$h$ be a singular metric with 
semi-positive curvature on $F$. 
Then for a $($non-zero$)$ section $s$ of a positive multiple $F^{m}$ 
satisfying $\sup_{X} |s|_{h^{m}} < \infty$, 
the multiplication map 
\begin{equation*}
\Phi_{s}: H^{q}(X, K_{X} \otimes F \otimes \I{h}) 
\xrightarrow{\otimes s} 
H^{q}(X, K_{X} \otimes F^{m+1} \otimes \I{h^{m+1}})
\end{equation*}
is $($well-defined and$)$ injective for any $q$. 
Here $\I{h}$ denotes the multiplier ideal sheaf 
associated to the singular metric $h$.  
\end{theo}

\begin{rem}
The multiplication map is well-defined 
thanks to the assumption 
of $\sup_{X} |s|_{h^{m}} < \infty $. 
When $h$ is a metric with minimal singularities 
on $F$, this assumption is automatically satisfied 
for any section $s$ of $F^{m}$ 
(see \cite{Dem} for the concept of metrics with minimal singularities).   
\end{rem}

When we consider the problem of extending 
(holomorphic) sections from subvarieties to the ambient space, 
we need to refine the above formulation (see Theorem \ref{main-inj2}).  
Our injectivity theorem 
can be seen as an improvement of 
\cite{Eno90}, \cite{Fuj12a}, \cite{Kol86}, \cite{Mat14}. 
For the proof, 
we take an analytic approach 
for the cohomology groups with coefficients in 
$K_{X} \otimes F \otimes \I{h}$, 
which includes techniques of 
\cite{Eno90}, \cite{Fuj12a}, \cite{Mat13a}, \cite{Mat14}, 
\cite{Ohs04}, \cite{Tak97}. 
The proof is based on a technical combination of 
the $L^{2}$-method for the $\dbar$-equation 
and the theory of harmonic integrals. 
To handle transcendental (non-algebraic) singularities, 
after regularizing a given singular metric, we 
investigate the asymptotic behavior of the harmonic forms 
with respect to a family of the regularized metrics. 
Moreover we establish a method to obtain  
$L^{2}$-estimates of solutions of the $\dbar$-equation  
by using the $\rm{\check{C}}$ech complex. 
See subsection \ref{S-Inj} for more details.

\subsection{Applications to the vanishing theorem}

Theorem \ref{main-inj} is  formulated for 
singular metrics with transcendental (non-algebraic) singularities, 
which is one of the advantages of our injectivity theorem. 
For example, 
metrics with minimal singularities are important objects, 
but they do not always have algebraic singularities. 
By applying Theorem \ref{main-inj} to them, 
we can obtain an injectivity theorem 
for nef and abundant line bundles (Corollary \ref{good}) 
and Nadel type vanishing theorems  
(Theorem \ref{gen} and Corollary \ref{main-co}).

It is natural to expect 
the same conclusion as in Theorem \ref{Kol} 
under the weaker assumption that $F$ is nef. 
However there is a counterexample  
to the injectivity theorem for nef line bundles 
(see for example \cite[Example 5.1]{Fuj12b}). 
If $F$ is nef and abundant (that is, the numerical dimension  
agrees with the Kodaira dimension), 
the line bundle $F$ 
admits a metric $h_{\min}$ with minimal singularities satisfying  $\I{h_{\min}^{m}}= \mathcal{O}_{X}$ for any $m>0$. 
This follows from \cite[Proposition 2.1]{Kaw85}. 
Hence  
Theorem \ref{main-inj} leads to the following corollary.

\begin{cor}[\text{\cite[Corollary 1.5]{Mat13b}}]\label{good}
Let $F$ be a nef and abundant line bundle 
on a compact K\"ahler manifold $X$. 
Then for a $($non-zero$)$ section $s$ of a positive multiple 
$F^{m}$ of the line bundle $F$, 
the multiplication map 
induced by the tensor product with $s$ 
\begin{equation*}
\Phi_{s}: H^{q}(X, K_{X} \otimes F) 
\xrightarrow{\otimes s} 
H^{q}(X, K_{X} \otimes F^{m+1} )
\end{equation*}
is injective for any $q$. 
\end{cor}

The same statement was proved in \cite{Fuj12a}, 
and a similar conclusion was proved in 
\cite{EP08}, \cite{EV92} by different methods 
when $X$ is a projective variety.
It is worth pointing out that Theorem \ref{Kol} is 
not sufficient to obtain Corollary \ref{good}. 
This is because   
the above metric $h_{\min}$ may not be smooth and 
not have algebraic singularities 
even if $F$ is nef and abundant  
(see for example \cite[Example 5.2]{Fuj12b}).

As another application of Theorem \ref{main-inj}, 
we obtain a Nadel type vanishing theorem 
(Theorem \ref{gen}) and its corollary (Corollary \ref{main-co}).

\begin{theo}[\text{\cite[Theorem 3.21]{Mat13b} 
cf. \cite[Theorem 5.2]{Mat14}}]
\label{gen}
Let $F$ be a line bundle 
on a smooth projective variety $X$ and 
$h$ be a singular metric with semi-positive curvature on $F$. 
Then 
\begin{equation*}
H^{q}(X, K_{X}\otimes F \otimes \I{h}) = 0
\hspace{0.4cm} {\text{for}}\ {\text{any}}\ 
q > \dim X-\kappa_{\rm{bdd}}(F, h).
\end{equation*}
See  subsection \ref{proofgen}
or \cite[Definition 5.1]{Mat14} for the definition of 
the bounded Kodaira dimension $\kappa_{\rm{bdd}}(F, h)$. 
\end{theo}

\begin{cor}[\text{\cite[Corollary 1.6]{Mat13b} 
cf. \cite[Theorem 1.2]{Mat13a}}]
\label{main-co}
Let $F$ be a line bundle on a smooth projective 
variety $X$ and  
$h_{\min}$ be 
a singular metric with minimal singularities on $F$.
Then  
\begin{equation*}
H^{q}(X, K_{X} \otimes F \otimes \I{h_{\min}}) = 0 
\hspace{0.4cm} {\text{for}}\ {\text{any}}\ 
q > \dim X-\kappa(F).
\end{equation*}
Here $\kappa(F)$ denotes the Kodaira dimension of $F$. 
\end{cor}

Since multiplier ideal sheaves are coherent ideal sheaves, 
the family of multiplier ideal sheaves 
$\{ \mathcal{I}(h^{1+\delta})\}_{\delta>0}$ has 
the maximal element, which we denote by $\mathcal{I}_{+}(h)$ 
(see \cite{DEL00} for more details). 
In \cite{Cao12}, Cao gave 
striking results on 
the Nadel vanishing theorem for the 
cohomology groups with coefficients in 
$K_{X}\otimes F \otimes \mathcal{I}_{+}(h)$.     
However, the Nadel vanishing theorem 
for $K_{X}\otimes F \otimes \I{h_{\min}}$ 
is non-trivial even if $F$ is big. 
In fact, it was first proved in \cite{Mat13a} 
when $F$ is big.

It is of interest to ask whether or not 
$\mathcal{I}_{+}(\varphi)$ agrees with $\mathcal{I}(\varphi)$ 
for a plurisubharmonic (psh for short) function $\varphi$,  
which was first posed in \cite{DEL00}. 
We can easily see that 
$\mathcal{I}_{+}(\varphi) = \I{\varphi}$ holds 
if $\varphi$ has algebraic singularities,  
but $h_{\min}$ unfortunately does not always have 
algebraic singularities.  
It is a natural problem  
related to the (strong) openness conjecture of Demailly-Koll\'ar (see \cite{DK01}), but it had been an open problem.  
Recently, Guan-Zhou affirmatively solved 
the openness conjecture in \cite{GZ13}. 
Although their celebrated results imply Theorem \ref{gen}, 
we believe that our techniques are still of interest, 
since they bring a quite different viewpoint and 
have further applications.

\subsection{Applications to the extension theorem}
In this subsection, we give an extension theorem for  
pluri-canonical sections of log pairs. 
Our motivation is the abundance conjecture, 
which is  one of 
the most important conjectures in the classification theory of algebraic varieties. 
From now on, we freely use the standard notation in 
\cite{BCHM10}, \cite{kamama}, \cite{KM} 
and  further we interchangeably use the words \lq \lq Cartier divisors", \lq \lq line bundles", \lq \lq invertible sheaves".

\begin{conj}[Generalized abundance conjecture]\label{abun conj}Let $X$ be a normal projective variety and $\Delta$ be an effective $\mathbb{Q}$-divisor such that $(X,\Delta)$ is a klt pair. 
Then  
$\kappa(K_X+\Delta)=\kappa_{\sigma}(K_X+\Delta)$.   In particular, if $K_X+\Delta$ is nef, then it is semi-ample. 
See  \cite{Nak} for the definition of $\kappa(\cdot)$ and $\kappa_{\sigma}(\cdot)$.
\end{conj}

Toward the abundance conjecture, 
we need to study the non-vanishing conjecture  
and the extension conjecture (see \cite{DHP13}, \cite[Introduction]{Fuj00}, \cite[Section 5]{FG14}).  
We study 
the following extension conjecture for dlt pairs 
formulated in \cite[Conjecture 1.3]{DHP13}:

\begin{conj}[{Extension conjecture for dlt pairs}]\label{c-dlt} 
Let $X$ be a normal projective variety 
and $S + B$ be an effective $\mathbb Q$-divisor
satisfying the following assumptions\,$:$ 
\begin{itemize}
\item $(X, S+B)$ is a dlt pair.
\item $\lfloor S+B \rfloor =S$.  
\item $K_X+S+B$ is nef. 
\item $K_X+S+B$ is $\mathbb Q$-linearly equivalent to 
an effective divisor $D$ with $S \subseteq \Supp D \subseteq \Supp\,(S+B)$. 
\end{itemize}
Then the restriction map 
$$H^0(X,\mathcal O _X(m(K_X+S+B)))\to H^0(S,  \mathcal O _S(m(K_X+S+B)))  $$
is surjective for sufficiently divisible integers $m\geq 2$. 
\end{conj}

When $S$ is a normal irreducible variety 
(that is, 
$(X, S+B)$ is a plt pair), 
Demailly-Hacon-P\u{a}un proved the above conjecture in \cite{DHP13}
by using technical methods based on a version of the Ohsawa-Takegoshi $L^{2}$-extension theorem.
This result can be seen as an extension theorem 
for plt pairs.

By applying 
Theorem \ref{main-inj2} instead of the Ohsawa-Takegoshi theorem 
to the extension problem,  
we prove the following extension theorem for {\textit{dlt pairs}}. 
Thanks to the injectivity theorem,  
we can obtain some extension theorems for not only plt paris but also dlt pairs. 
This is an advantage of our approach. 
Even if $K_{X}+\Delta$ is semi-positive 
(that is, it admits a smooth metric with semi-positive curvature), 
it seems to be very impossible to prove the extension theorem for dlt pairs by the Ohsawa-Takegoshi theorem at least in the current situation, 
and thus we need our injectivity theorem 
(Theorem \ref{main-inj2}).

\begin{theo}[\text{\cite[Corollary 4.5]{GM13}}]
\label{intro-ext}
Let $X$ be a compact K\"ahler manifold and 
$S+B$ be an effective $\mathbb Q$-divisor with the 
following assumptions\,$:$ 
\begin{itemize}
\item $S+B$ is a simple normal crossing divisor with 
$0 \leq S+B \leq 1$ and $\lfloor S+B \rfloor =S$.
\item $K_{X}+S+B$ is ${\mathbb{Q}}$-linearly 
equivalent to an effective divisor $D$ 
with $S \subseteq\mathrm{Supp}\, D$. 
\item $K_{X}+S+B$ admits a singular metric $h$ with semi-positive curvature. 
\item  The Lelong number $\nu(h,x)$ is equal to $0$ at every point 
$x \in S$. 
\end{itemize}
Then, for an integer $m \geq 2$ with Cartier divisor $m(K_{X}+S+B)$,  
every section $u \in H^{0}(S, \mathcal{O}_{S}(m(K_{X}+S+B)))$ 
can be extended to a section in $H^{0}(X, \mathcal{O}_{X}(m(K_{X}+S+B)) )$. 
\end{theo}

In particular, Conjecture \ref{c-dlt} is affirmatively solved  
under the assumption that $K_{X}+ \Delta$ admits a singular metric 
whose curvature is semi-positive and Lelong number is identically zero. 
This assumption is stronger than the assumption that $K_{X}+\Delta$ is nef, but 
weaker than the assumption that $K_{X}+\Delta$ is semi-positive. 
Let us observe that Verbitsky proved the non-vanishing conjecture 
on hyperK\"ahler manifolds (holomorphic symplectic manifolds)  
under the same assumption (see \cite{Ver10}).

As compared to Conjecture \ref{c-dlt}, one of our advances  
has been to remove the condition  $\Supp\,D\subseteq \Supp (S+B) $. 
As a benefit of removing the condition  $\Supp\,D\subseteq \Supp (S+B) $ in Conjecture \ref{c-dlt}, we can run the minimal model program 
while preserving a good condition for metrics 
(cf. \cite[Section 8]{DHP13}, \cite[Theorem 5.9]{FG14}). 
By applying the above theorem and techniques of the  minimal model program, 
we obtain results related to 
the abundance conjecture 
(see \cite{GM13} for more details).

\subsection*{Acknowledgement}

The author obtained an opportunity of discussion on 
the injectivity theorem and extension problem 
when he attended 
the conference \lq \lq The 10th Korean Conference in Several Complex Variables". 
He is grateful to the organizers. 
He would also like to thank the referee for
carefully reading the paper and for suggestions.
He is partially supported by the Grant-in-Aid for 
Young Scientists (B) $\sharp$25800051 from JSPS.

\section{Proof of the Main Results}\label{S-Proof}

\subsection{Proof of Theorem \ref{main-inj2}}\label{S-Inj}
In this subsection, 
we give a proof of the following theorem, 
which is an improvement of Theorem \ref{main-inj} 
to obtain Theorem \ref{intro-ext}.

\begin{theo}\label{main-inj2}
Let $(F, h_{F})$ and $(L, h_{L})$ 
be $($singular$)$ hermitian line bundles 
with semi-positive curvature 
on a compact K\"ahler manifold $X$. 
Assume that there exists 
an effective $\mathbb{R}$-divisor $\Delta$ with 
\begin{align*}
h_{F}=h_{L}^{a}\cdot h_{\Delta},  
\end{align*}
where $a$ is a positive real number and  $h_{\Delta}$  is 
the singular metric defined by $\Delta$.

Then for a $($non-zero$)$ section $s$ of $L$
satisfying $\sup_{X} |s|_{h_{L}} < \infty$, 
the multiplication map 
\begin{equation*}
\Phi_{s}: H^{q}(X, K_{X} \otimes F \otimes \I{h_{F}}) 
\xrightarrow{\otimes s} 
H^{q}(X, K_{X} \otimes F\otimes L \otimes \I{h_{F} h_{L}})
\end{equation*}
is $($well-defined and$)$ injective for any $q$.  
\end{theo}

\begin{rem}\label{rem-inj}

(1) The case of $\Delta=0$  
corresponds to Theorem \ref{main-inj}. 
\vspace{0.1cm}\\
$(2)$ If $h_{L}$ and $h_{F}$ are smooth on a Zariski open set, 
the same conclusion holds under the weaker assumption of 
$\sqrt{-1}\Theta_{h_{F}}(F) \geq a \sqrt{-1}\Theta_{h_{L}}(L)$ 
$($see \cite[Theorem 1.5]{Mat14}$)$. 
\end{rem}

\begin{proof}
We give here only the strategy of the proof.  
See \cite{Mat13b}, \cite{GM13} for the precise proof. 
First of all, we recall Enoki's method to generalize  
Koll\'ar's injectivity theorem, which gives a proof of   
the special case where $h_{L}$ is smooth and $\Delta=0$. 
In this case, 
the cohomology group $H^{q}(X, K_{X} \otimes F)$ 
is isomorphic to the space 
of the harmonic forms with respect to $h_{F}$
\begin{equation*}
\mathcal{H}^{n, q}(F)_{h_{F}}:= \{u \mid u 
\text{ is a smooth $F$-valued $(n,q)$-form on } X 
\text{ such that } \dbar u  = D^{''*}_{h_{F}} u =0.  \}, 
\end{equation*}
where  $D^{''*}_{h_{F}}$ is the adjoint operator of 
the $\dbar$-operator. 
For an arbitrary harmonic form 
$u \in \mathcal{H}^{n, q}(F)_{h_{F}}$,  
we can conclude that $D^{''*}_{h_{F}h_{L}} su =0$  
from the semi-positivity of the curvature 
and $h_{F}=h^{a}_{L}$.
This step heavily depends on the semi-positivity 
of the curvature.  
This implies that the multiplication map $\Phi_{s}$ 
induces the map from $\mathcal{H}^{n, q}(F)_{h_{F}}$ to 
$\mathcal{H}^{n, q}(F\otimes L)_{h_{F}h_{L}}$, and thus 
the injectivity is obvious.

When $h_{L}$ is smooth on a Zariski open set, 
the cohomology group 
$H^{q}(X, K_{X} \otimes F \otimes \I{h})$ 
is isomorphic to the space of harmonic forms on 
the Zariski open set. 
Therefore we can give a proof similar to Enoki's proof
thanks to the semi-positivity of the curvature 
(see \cite[Theorem 1.5]{Mat14}).

In our situation, we must consider   
singular metrics with transcendental 
(non-algebraic) singularities. 
It is quite difficult 
to directly handle transcendental singularities, 
and thus, in Step 1, we approximate a given singular metric $h_{F}$ by metrics 
$\{h_{\e} \}_{\e>0}$ that are smooth on a Zariski open set. 
Then we represent a given cohomology class in 
$H^{q}(X, K_{X} \otimes F \otimes \I{h_{F}})$ 
by the associated harmonic form  
$u_{\e}$ with respect to $h_{\e}$ on the Zariski open set. 
We want to show that $s u_{\e}$ is also harmonic 
by using the same method as Enoki. 
However, the same argument as in \cite{Eno90} fails  
since the curvature of $h_{\e}$ is not semi-positive. 
For this reason, in Step2, 
we investigate the asymptotic behavior of 
the harmonic forms $u_{\e}$ with respect to  
a family of the regularized metrics 
$\{ h_{\e} \}_{\e>0}$.  
Then we show that the $L^{2}$-norm 
$\| D^{''*}_{ h_{\e} h_{L, \e} }su_{\e}\|$ converges to zero 
as $\e$ tends to zero, 
where $h_{L, \e}$ is a suitable approximation of $h_{L}$. 
Further, in Step 3, we construct solutions $\gamma_{\e}$ of 
the $\dbar$-equation $\dbar \gamma_{\e} = su_{\e}$ 
such that   
the $L^{2}$-norm $\| \gamma_{\e} \|$  
is uniformly bounded,  
by applying the $\rm{\check{C}}$ech complex with 
the topology induced by the local $L^{2}$-norms. 
In Step 4, we see that 
\begin{equation*}
\| su_{\e} \| ^{2} = 
\lla su_{\e}, \dbar\gamma_{\e} \rra
\leq \| D^{''*}_{h_{\e} h_{L, \e}} su_{\e}\|
\| \gamma_{\e} \| \to 0 \quad {\text{as }} \e \to 0.
\end{equation*}
From these observations, we conclude  that 
$u_{\e}$ converges to zero in a suitable sense. 
This completes the proof. 
\vspace{0.2cm} \\
{\bf{Step 1 (The equisingular approximation of $h_{F}$)}}
\vspace{0.1cm} \\
Throughout the proof, we fix a K\"ahler form $\omega$ on $X$. 
For the proof, 
we want to apply the theory of harmonic integrals, 
but the metric $h_{F}$ may not be smooth. 
For this reason, we approximate $h_{F}$ by metrics 
$\{ h_{\e} \}_{\e>0}$ that are smooth on a Zariski open set. 
By \cite[Theorem 2.3]{DPS01}, we can obtain  
metrics $\{ h_{\e} \}_{\e>0}$ on $F$ satisfying the following properties:
\begin{itemize}
\item[(a)] $h_{\e}$ is smooth on $Y:=X \setminus Z$, 
where $Z$ is a subvariety independent of $\e$.
\item[(b)]$h_{\e_{2}} \leq h_{\e_{1}} \leq h_{F}$ 
holds for any 
$0< \e_{1} < \e_{2} $.
\item[(c)]$\I{h_{F}}= \I{h_{\e}}$.
\item[(d)]$\sqrt{-1} \Theta_{h_{\e}}(F) \geq -\e \omega$. 
\end{itemize}
See \cite[Theorem 2.3]{Mat13b} for property (a). 
By \cite[Lemma 3.1]{Fuj12a}, 
we obtain a K\"ahler form $\ome$ on $Y$   
satisfying the following properties: 
\begin{itemize}
\item[(A)] $\ome$ is a complete K\"ahler form on $Y$.
\item[(B)] There exists a bounded function $\Psi$ such that 
$\ome = \ddbar \Psi$ on a neighborhood of $z\in Z$. 
\item[(C)] $\ome \geq \omega $.
\end{itemize}

In the proof, 
we mainly consider harmonic forms on $Y$ 
with respect to $h_{\e}$ and $\ome$. 
Let $L_{(2)}^{n, q}(Y, F)_{h_{\e}, \ome}$ be   
the space of $L^{2}$-integrable $F$-valued $(n,q)$-forms 
$\alpha$ with respect to the inner product $\|\cdot \|_{h_\e, \ome}$
defined by 
\begin{equation*}
\|\alpha \|^{2}_{h_\e, \ome}:= \int_{Y} 
|\alpha |^{2}_{h_{\e}, \ome}\ \ome^{n}. 
\end{equation*}
Then we have the following orthogonal decomposition: 
\begin{equation*}
L_{(2)}^{n, q}(Y, F)_{h_{\e}, \ome}
=
{\rm{Im}}\,\dbar
\oplus
\mathcal{H}^{n, q}(F)_{h_{\e}, \ome}
\oplus {\rm{Im}}\,D^{''*}_{h_{\e}}.   
\end{equation*}
Here the operator $D^{'*}_{h_{\e}}$ 
(resp. $D^{''*}_{h_{\e}}$) denotes   
the closed extension of the formal adjoint of the 
$(1,0)$-part $D^{'}_{h_{\e}}$ (resp. $(0,1)$-part $D^{''}_{h_{\e}}=\dbar$) 
of the Chern connection $D_{h_{\e}}=D^{'}_{h_{\e}}+ D^{''}_{h_{\e}}$.  
Further $\mathcal{H}^{n, q}(F)_{h_{\e}, \ome}$ denotes   
the space of harmonic forms with respect to 
$h_{\e}$ and $\ome$, namely 
\begin{equation*}
\mathcal{H}^{n, q}(F)_{h_{\e}, \ome}:= 
\{\alpha   \mid \alpha 
\text{ is an } F\text{-valued } (n,q)\text{-form with  }
\dbar \alpha= D^{''*}_{h_{\e}}\alpha=0.    \}. 
\end{equation*}
A harmonic form in $ \mathcal{H}^{n, q}(F)_{h_{\e}, \ome}$ 
is smooth by the regularity  theorem for elliptic operators. 
These results are known to specialists.  
The precise proof of them   
can be found in \cite[Claim 1]{Fuj12a}.

Take an arbitrary cohomology class
$\{u \} \in H^{q}(X, K_{X} \otimes F \otimes \I{h_{F}})$ 
represented by an $F$-valued 
$(n, q)$-form $u$ with $\|u \|_{h_{F}, \omega} < \infty$. 
In order to prove that 
the multiplication map $\Phi_{s}$ is injective, 
we assume that the cohomology class of $su$ 
is zero in 
$H^{q}(X, K_{X}\otimes F\otimes L \otimes \I{h_{F}h_{L}})$. 
Our goal is to show that 
the cohomology class of $u$ is actually zero 
under this assumption.

By the inequality  $\|u\|_{h_{\e}, \ome} \leq 
\|u \|_{h_{F}, \omega}< \infty$,   
we can obtain 
$u_{\e} \in \mathcal{H}^{n, q}(F)_{h_{\e}, \ome}$ and 
$v_{\e} \in L_{(2)}^{n,q-1}(Y, F)_{h_{\e}, \ome}$ such that 
\begin{equation*}
u=u_{\e}+\dbar v_{\e}. 
\end{equation*}
Note that  
the component of ${\rm{Im}} D^{''*}_{h_{\e}}$ is zero 
since $u$ is $\dbar$-closed.

At the end of this step, we explain  
the strategy of the proof. 
In Step 2, we show that 
$\|D^{''*}_{h_{\e} h_{L, \e}} 
s u_{\e} \|_{ h_{\e} h_{L, \e}, \ome}$ 
converges to zero 
as $\e$ tends to zero. 
Here $h_{L, \e}$ is the singular metric on $L$ defined by 
\begin{align*}
h_{L, \e}:= h_{\e}^{1/a}\, h_{\Delta}^{-1/a}. 
\end{align*}
Since the cohomology class of $su$ is zero, 
there are solutions $\gamma_{\e}$ of the $\dbar$-equation 
$\dbar \gamma_{\e} = s u_{\e}$. 
For the proof, we need to obtain $L^{2}$-estimates of them. 
In Step 3, 
we construct solutions $\gamma_{\e}$ of the $\dbar$-equation   
$\dbar \gamma_{\e} = s u_{\e}$ such that  
the norm $\| \gamma_{\e} \|_{h_{\e} h_{L, \e}, \ome}$ 
is uniformly bounded. 
Then we have 
\begin{equation*}
\|su_{\e} \|^{2}_{h_{\e} h_{L, \e}, \ome} \leq 
\|D^{''*}_{h_{\e} h_{L, \e}} s u_{\e} \|_{h_{\e} h_{L, \e}, \ome}
\| \gamma_{\e} \|_{h_{\e} h_{L, \e}, \ome}.
\end{equation*} 
By Step 2 and Step 3, we can conclude that  
the right hand side goes to zero as 
$\e$ tends to zero. 
In Step 4, from this convergence, we prove that 
$u_{\e}$ converges to zero  
in a suitable sense, 
which implies that the cohomology class of $u$ is zero. 
\vspace{0.1cm}\\
\vspace{0.2cm} \\
{\bf{Step 2 (A generalization of Enoki's proof)}}
\vspace{0.1cm}\\
By generalizing Enoki's method, 
in Step 2, we prove the following proposition:  
\begin{prop} \label{D''}
As $\e$ tends to zero, the norm 
$\|D^{''*}_{h_{\e} h_{L, \e}} s u_{\e} \|_{h_{\e} h_{L, \e}, \ome}$ converges  
to zero. 
\end{prop}

The same argument as in \cite{Eno90} fails 
since the curvature of $h_{\e}$ is not semi-positive, 
and further property (d) is not sufficient 
for the proof of the proposition 
since there is counterexample to the injectivity theorem 
for nef line bundles. 
To overcome these difficulties,  
we first see 
the following inequality: 
\begin{equation}\label{ine2}
\|u_{\e} \|_{h_{\e}, \ome} 
\leq \|u \|_{h_{\e}, \ome} 
\leq \|u \|_{h, \omega}.  
\end{equation}
This inequality and properties (b), (c) imply the proposition. 
This step can be considered as a generalization of Enoki's method. 
\vspace{0.2cm} \\
{\bf{Step 3 (A construction of 
solutions of the $\dbar$-equation via 
the $\bf\rm{\check{C}}$ech complex)}}
\vspace{0.1cm} \\
In Step 3, we construct solutions of the $\dbar$-equation 
with suitable $L^{2}$-norm   
by using the $\bf\rm{\check{C}}$ech complex. 
\begin{prop}\label{sol-1}
There exist $F$-valued $(n, q-1)$-forms $\alpha_{\e}$ 
on $Y$ 
satisfying the following properties\,$:$ 
\begin{equation*}
{\rm{(1)}}\hspace{0.2cm} \dbar \alpha_{\e}=u-u_{\e}.   
\quad 
{\rm{(2)}}\hspace{0.2cm} \text{The norm }
\|\alpha_{\e} \|_{h_{\e}, \ome} 
 \text{ is uniformly bounded}. 
\end{equation*}
\end{prop}
\begin{rem}
We have already known that there exist solutions $\alpha_{\e}$ of 
the $\dbar$-equation $\dbar \alpha_{\e}=u-u_{\e} $ 
since $u-u_{\e} \in {\rm{Im}} \dbar$. 
However, for the proof of the main theorem, 
we need to construct solutions with uniformly bounded $L^{2}$-norm. 
\end{rem}
The strategy of the proof is as follows:  
The main idea of the proof is to convert the $\dbar$-equation 
$\dbar \alpha_{\e}=u-u_{\e}$
to the equation $\delta V_{\e} = S_{\e}$ of 
the coboundary operator $\delta$
in the space of cochains 
$C^{\bullet}(K_{X}\otimes F \otimes \I{h_{\e}})$, 
by using the $\rm{\check{C}}$ech complex and pursuing  
the De Rham-Weil isomorphism. 
Here the $q$-cochain $S_{\e}$ is constructed from $u-u_{\e}$.  
In this construction, we locally solve the $\dbar$-equation. 
The important point is that the space 
$C^{\bullet}(K_{X}\otimes F \otimes \I{h_{\e}})$ is independent of 
$\e$ thanks to property (c) of $h_{\e}$
although the $L^{2}$-space 
$L^{n,q}_{(2)}(Y, F)_{h_{\e}, \ome}$ depends on $\e$. 
Since $\|u-u_{\e}\|_{h_{\e}, \ome}$ is uniformly bounded, 
we can observe that $S_{\e}$ converges to 
some $q$-coboundary in $C^{q}(K_{X}\otimes F \otimes \I{h})$ 
with the topology induced by 
the local $L^{2}$-norms with respect to $h$.
Further we can observe that 
the coboundary operator $\delta$ is an open map.
Then by these observations we construct solutions $V_{\e}$ of 
the equation $\delta V_{\e} = S_{\e}$ 
with uniformly bounded norm. 
Finally, by using a partition of unity, 
we conversely 
construct $\alpha_{\e} \in L^{n,q-1}_{(2)}(Y, F)_{h_{\e}, \ome}$
from $S_{\e}$ satisfying the properties in Proposition \ref{sol-1}. 
This proof gives a new method to 
obtain $L^{2}$-estimates of solutions of 
the $\dbar$-equation.

\hspace{0.1cm}\\ \ 
{\bf{Step 4 (The limit of the harmonic forms)}}
\vspace{0.1cm}\\ \ 
In Step 4, we investigate the limit of $u_{\e}$ and 
complete the proof. 
By Step 2 and Step 3, we have  
\begin{equation*}
\|su_{\e} \|^{2}_{h_{\e} h_{L, \e}, \ome} \leq 
\|D^{''*}_{h_{\e} h_{L, \e}} s u_{\e} \|_{h_{\e} h_{L, \e}, \ome}
\| \gamma_{\e} \|_{h_{\e} h_{L, \e}, \ome} \to 0 
\quad  \text{as} \quad  \e   \to 0.
\end{equation*} 
From this convergence, we can show that 
$u_{\e}$ converges to zero  in a suitable sense, 
which implies that the cohomology class $\{u \}$ of $u$ is zero
in $H^{q}(X, K_{X} \otimes F \otimes \I{h_{\e}})$. 
By property (c), we obtain the conclusion of Theorem \ref{main-inj2}. 
\end{proof}

\subsection{Proof of Theorem \ref{gen}}
\label{proofgen}
In this subsection, we give 
a proof of Theorem \ref{gen} by using Theorem \ref{main-inj} and 
\cite[Theorem 4.1]{Mat14}.

\vspace{0.2cm}
\hspace{-0.6cm}
\textit{Proof of Theorem \ref{gen}.$)$}
We consider 
the space of sections with bounded norm 
defined by  
\begin{equation*}
H_{{\rm{bdd}}, h^m}^{0}(X, F^m): =
\{ s \in H(X, F^m)\ \big| \ \sup_{X} |s|_{h^m} < \infty \}.
\end{equation*}
The {\textit{generalized Kodaira dimension}} 
$\kappa_{\rm{bdd}}(F, h)$ of $(F,h)$ is defined to be $-\infty$ 
if $H_{{\rm{bdd}}, h^{m}}^{0}(X,F^{m})=0$ for any 
$m>0$. 
Otherwise, $\kappa_{\rm{bdd}}(F, h)$ is defined 
by 
\begin{equation*}
\kappa_{\rm{bdd}}(F, h): = 
\sup \{ k \in \mathbb{Z}\mid \limsup_{m \to \infty} 
\dim H_{{\rm{bdd}}, h^{m}}^{0}(X, F^{m})\big/ m^{k} > 0 \}. 
\end{equation*}
For a contradiction, 
we assume that there exists a non-zero cohomology class 
$\alpha \in H^{q}(X, K_{X}\otimes F \otimes \I{h})$. 
If sections  
$\{s_{i}\}_{i=1}^{N}$ in 
$H_{{\rm{bdd}}, h^{m}}^{0}(X, F^{m})$ are 
linearly independent, 
then  $\{s_{i} \alpha \}_{i=1}^{N}$ 
are also  linearly independent in 
$ H^{q}(X, K_{X}\otimes F^{m+1} \otimes \I{h^{m+1}})$. 
Indeed, if 
$\sum_{i=1}^{N} c_{i} s_{i} \alpha =0$  
for some $c_{i} \in \mathbb{C}$, 
then we know $\sum_{i=1}^{N} c_{i} s_{i}=0$ 
by Theorem \ref{main-inj}. 
Since $\{s_{i}\}_{i=1}^{N}$ are linearly independent, 
we have $c_{i}=0$ for any $i=1,2,\dots N$. 
This yields   
\begin{equation*}
\dim H_{{\rm{bdd}}, h^{m}}^{0}(X, F^{m}) \leq 
\dim H^{q}(X, K_{X}\otimes F^{m+1} \otimes \I{h^{m+1}}). 
\end{equation*}
On the other hand, 
by \cite[Theorem 4.1]{Mat14}, 
we have  
\begin{equation*}
\dim H^{q}(X, K_{X}\otimes F^{m} \otimes \I{h^{m}})= O(m^{\dim X-q}) 
\quad  \text{as} \quad  m   \to \infty
\end{equation*}
for any $q\geq 0$ 
(cf. \cite[(6.18) Lemma]{Dem}). 
If $q > \dim X-\kappa_{\rm{bdd}}(F, h)$, this is a contradiction.  
\qed

\subsection{Proof of Theorem \ref{intro-ext}}
In this subsection, we give a proof of Theorem \ref{intro-ext}.

\vspace{0.2cm}
\hspace{-0.6cm}
\textit{Proof of Theorem \ref{intro-ext}.$)$}  
For simplicity, we put $\Delta:= S+B$ and $G:=m(K_{X}+ \Delta)$. 
We may assume 
the additional assumption of  
$h \leq h_{D}$, 
where $h_{D}$ is the singular metric on $\mathcal{O}_{X}(K_{X}+\Delta)$ 
defined by the effective divisor $D$. 
Indeed, for a smooth metric $g$ on $\mathcal{O}_{X}(K_{X}+\Delta)$ and 
an $L^{1}$-function $\varphi$ (resp. $\varphi_{D}$)  
with $h= g\, e^{-\varphi} $ (resp. $h_{D}= g\, e^{-\varphi_{D}} $),  
the metric defined by $g\, e^{-\max( \varphi, \varphi_{D})}$
satisfies the assumptions again.

Consider the following exact sequence: 
\begin{equation*}
0 \rightarrow \mathcal{O}_{X}(G - S) \otimes I (h^{m-1}h_{B}) 
\rightarrow \mathcal{O}_{X}(G ) \otimes \I{h^{m-1}h_{B}}
\rightarrow \mathcal{O}_{S}(G ) \otimes \I{h^{m-1}h_{B}} 
\rightarrow0.
\end{equation*}
We first prove 
the induced homomorphism
\begin{equation*}
H^{q}(X, \mathcal{O}_{X}(G - S) \otimes I (h^{m-1}h_{B})) 
\to 
H^{q}(X, \mathcal{O}_{X}(G) \otimes I (h^{m-1}h_{B})) 
\end{equation*}
is injective by our injectivity theorem. 
By the assumption on the support of $D$, we can take 
an integer $a>0$ such that $aD$ is a 
Cartier divisor and $S \leq a D$. 
Then we have the following commutative diagram:

\[\xymatrix{
& \hspace{-3.5cm} H^{q}(X, \mathcal{O}_{X}(G) \otimes I (h^{m-1}h_{B})) \supseteq \mathrm{Im}\,(+S)   
 \ar[d]^{+ (aD -S)}   
& \\
H^{q}(X, \mathcal{O}_{X}(G - S) \otimes I (h^{m-1}h_{B}))  
\ar[ru]^{+ S}   
\ar[r]_{\hspace{-1.0cm}+ aD}  & 
H^{q}(X, \mathcal{O}_{X}(G - S + aD) \otimes I (h^{a+m-1}h_{B})), & 
}\]
with a map $+S:H^{q}(X, \mathcal{O}_{X}(G - S) \otimes I (h^{m-1}h_{B}))  \to H^{q}(X, \mathcal{O}_{X}(G ) \otimes I (h^{m-1}h_{B})).$ 
In order to show 
that the upper map on right is injective,  
we prove that the horizontal map is injective 
as an application of Theorem \ref{main-inj2}.

By the definition of $G$, 
we have 
\begin{equation*}
G -S =m(K_{X}+\Delta)-S=K_{X} + (m-1)(K_{X}+ \Delta) + B. 
\end{equation*}
Then the line bundle $F:= \mathcal{O}_{X}((m-1)(K_{X}+ \Delta) + B)$ 
equipped with the metric $h_{F}:= h^{m-1} h_{B}$
and the line bundle $L:=\mathcal{O}_{X}(aD)$ 
equipped with the metric $h_{L}:= h^{{a}}$ 
satisfy the assumptions in Theorem \ref{main-inj2}. 
Indeed, we have $h_{F}=h_{L}^{(m-1)/a} h_{B}$
by the construction, and further 
the point-wise norm $|s_{aD}|_{h_{L}}$ is bounded on $X$ 
by the inequality $h \leq h_{D}$, 
where $s_{aD}$ is the natural section of $aD$.    
Therefore the horizontal map is injective 
by Theorem \ref{main-inj2}. 
By the assumption on the Lelong number of $h$, 
we can conclude that $\mathcal{O}_{S}\otimes \I{h^{m-1}h_{B}}=\mathcal{O}_{S}$. 
This follows from Skoda's lemma and H\"older's inequality.
This completes the proof. 
\qed

\section{Open Problems}\label{S-App}
In this section, we summarize and give open problems related to 
the topics mentioned in this survey. 

It is of interest to consider the injectivity theorem 
in the relative situation. 
The following problem is a relative version of Theorem \ref{main-inj}. 
For relative versions of the injectivity theorem and 
their applications, we refer the reader to \cite{Fuj12b}. 
In his paper \cite{Fuj12b}, 
Fujino affirmatively solved this problem 
under the assumption on the regularity of singular metrics, 
whose proof is based  on 
the Ohsawa-Takegoshi twisted version of Nakano's identity. 
To remove this assumption, 
it seems to be needed to use a combination of 
his method and the techniques of Theorem \ref{main-inj}.  

\begin{prob}[\text{cf. \cite[Problem 1.8]{Fuj12b}}]
\label{prob-rel}
Let $\pi:X \to Y$ be a surjective holomorphic map 
from K\"ahler manifold $X$ to a complex manifold $Y$, 
and $F$ be a line bundle on $X$ with 
a singular metric $h$ whose curvature is semi-positive. 
Then for a $($non-zero$)$ section $s$ of 
a positive multiple $F^{m}$ satisfying
$\sup_{X}|s|_{h^{m}}< \infty$, 
the multiplication map 
\begin{equation*}
\Phi_{s}: R^{q}\pi_{*}(K_{X} \otimes F \otimes \I{h}) 
\xrightarrow{\otimes s} 
 R^{q}\pi_{*}(K_{X} \otimes F^{m+1} \otimes \I{h^{m+1}}) 
\end{equation*}
is injective for any $q$. 
Here  $R^{q}\pi_{*}(\mathcal{F})$ denotes 
the higher direct image of a sheaf $\mathcal{F}$. 
\end{prob}

Theorem \ref{main-inj2} 
can be expected to hold under the weaker assumption made 
in the following problem. 
Indeed, this problem was affirmatively solved in \cite{Mat14} 
under s regularity assumption on singular metrics. 
It is also an interesting problem to consider the relative version 
of this problem in the same direction as Problem \ref{prob-rel}.

\begin{prob}[\text{cf. \cite[Theorem 1.5]{Mat14}, 
\cite[Theorem 1.2]{Fuj12a}}]
\label{prob-dif}
Let $(F, h_{F})$ and $(L, h_{L})$ 
be $($singular$)$ hermitian line bundles 
with semi-positive curvature 
on a compact K\"ahler manifold $X$. 
Assume there exists a positive real number $a$  
such that $\sqrt{-1}\Theta_{h_{F}}(F) \geq a \sqrt{-1}\Theta_{h_{L}}(L)$. 
Then the same conclusion as in Theorem \ref{main-inj2} holds. 
\end{prob}

Fujino proposed the following problem, 
which asks whether one can generalize 
the injectivity theorem for lc pairs proved by him. 
The main difficulty in studying this problem is 
that one must handle lc singularities by analytic methods.

\begin{prob}[cf. \text{\cite[Theorem 6.1]{Fuj11}}]
Let $D$ be a simple normal crossing divisor 
and $F$ be a semi-positive line bundle 
on a compact K\"ahler manifold $X$. 
Then, 
for a $($non-zero$)$ section $s$ of a positive multiple 
$F^{m}$ whose zero locus $s^{-1}(0)$ contains 
no lc centers of $(X,D)$, 
the multiplication map 
%induced by the tensor product with $s$ 
\begin{equation*}
\Phi_{s}: H^{q}(X, K_{X}\otimes F\otimes \mathcal{O}_{X}(D)) 
\xrightarrow{\otimes s} 
H^{q}(X, K_{X}\otimes F^{m+1}\otimes \mathcal{O}_{X}(D) )
\end{equation*}
is injective for any $q$. 
\end{prob}

For a nef line bundle $F$ on a smooth projective variety $X$, 
it can be proven that 
\begin{equation*}
\dim H^{q}(X, F^{m})= O(m^{\dim X -q})
\quad  \text{as} \quad  m   \to \infty.  
\end{equation*} 
When $X$ is merely supposed to be a compact K\"ahler manifold, 
the same conclusion can be expected. 
This was first posed by Demailly, 
and proved by Berndtsson under the stronger assumption 
that $F $ is semi-positive in \cite{Ber12}. 
The following problem was also proved in \cite{Mat14} 
when $X$ is a smooth projective variety.

\begin{prob}[\text{cf. \cite[Theorem 4.1]{Mat14}}]
Let $F$ be a line bundle on a compact K\"ahler manifold $X$ 
and $h$ be a singular metric with semi-positive curvature on $F$. 
Then, for any vector bundle $($or line bundle$)$ $M$, we have 
\begin{equation*}
\dim H^{q}(X, M\otimes F^{m}\otimes \I{h^{m}})= O(m^{\dim X -q})
\quad  \text{as} \quad  m   \to \infty. 
\end{equation*}
\end{prob}

%%%%%%%%%%%%%%%%%%%%%%%%%%%%%%%%%%%%%%%%%%%%%%%%%%%%%%%%%%%%%%%%%%%%%%%%%%%%%%%%%%%

%%%%%%%%%%%%%%%%%%%%%%%%%%%%%%%%%%%%%%%%%%%%%%%%%%%%%%%%%%%%%%%%%%%%%%%%%%%%%%%%%%%

\begin{thebibliography}{n}

\bibitem[Amb03]{Amb03}
F. Ambro.
{\textit{Quasi-log varieties.}} 
Tr. Mat. Inst. Steklova {\bf{240}} (2003), 
Biratsion. Geom. Linein. Sist. Konechno Porozhdennye Algebry, 220--239; 
translation in Proc. Steklov Inst. Math. (2003), no. 1 (240), 214--233. 


\bibitem[Amb12]{Amb12}
F. Ambro.   
\textit {An Injectivity Theorem.}
Compos. Math. {\bf{150}} (2014), no. 6, 999--1023.



\bibitem[BCHM10]{BCHM10}
C. Birkar, P. Cascini, C. D. Hacon, J. $\mathrm{M^{c}}$Kernan.
{\textit{Existence of minimal models for varieties of log general type.}} 
J. Amer. Math. Soc. {\bf{23}} (2010), 405--468.


\bibitem[Ber12]{Ber12}
B. Berndtsson. 
\textit{An eigenvalue estimate for the $\dbar$-Laplacian.} 
J. Differential Geom. {\bf {60}} (2002), no. 2, 295--313. 


\bibitem[Ber13]{Ber13}
B. Berndtsson.   
\textit{The openness conjecture for plurisubharmonic functions.}
Preprint, arXiv:1305.5781v1. 


\bibitem[Cao12]{Cao12}
J. Cao.  
\textit{Numerical dimension and a Kawamata-Viehweg-Nadel type vanishing theorem on compact K\"ahler manifolds.}
Preprint, arXiv:1210.5692v1, to appear in Composio. Math.



\bibitem[DEL00]{DEL00}
J.-P. Demailly, L. Ein, R. Lazarsfeld.  
\textit{A subadditivity property of multiplier ideals.}
Michigan Math. J. {\bf 48} (2000), 137--156.


\bibitem[Dem]{Dem}
J.-P. Demailly.   
\textit {Analytic methods in algebraic geometry.}
Surveys of Modern Mathematics {\bf{1}}, 
International Press, Somerville, 
MA; Higher Education Press, Beijing, (2012).


\bibitem[Dem-book]{Dem-book}
J.-P. Demailly.   
\textit {Complex analytic and differential geometry.}
Lecture Notes on the web page of the author.


\bibitem[Dem82]{Dem82}
J.-P. Demailly.  
\textit{Estimations $L^{2}$ pour l'op\'erateur $\overline{\partial}$ d'un 
fibr\'e vectoriel holomorphe semi-positif au-dessus d'une vari\'et\'e k\"ahl\'erienne compl\`ete.}
Ann. Sci. \'Ecole Norm. Sup(4). {\bf{15}} (1982), 457--511. 


\bibitem[DHP13]{DHP13}
J.-P. Demailly, C. D. Hacon, M. P\u{a}un.  
{\textit{Extension theorems, non-vanishing and the existence of good minimal models.}} 
Acta Math. {\bf 210} (2013), 203--259.


\bibitem[DK01]{DK01}
J.-P. Demailly, J. Koll\'ar.  
\textit{Semicontinuity of complex singularity exponents and K\"ahler-Einstein metrics on Fano orbifolds.}
Ann. Sci. \'Ecole Norm. Sup(4) {\bf {34}} (2001), 525--556.


\bibitem[DPS01]{DPS01}
J.-P. Demailly, T. Peternell, M. Schneider.  
\textit{Pseudo-effective line bundles on compact K\"ahler manifolds.}
International Journal of Math. {\bf{6}} (2001), 689--741. 


\bibitem[Eno90]{Eno90}
I. Enoki. 
\textit{Kawamata-Viehweg vanishing theorem for compact 
K\"ahler manifolds.}
Einstein metrics and Yang-Mills connections (Sanda, 1990), 
59--68. 


\bibitem[EP08]{EP08}
L. Ein, M. Popa. 
\textit{Global division of cohomology classes via injectivity.}
Special volume in honor of Melvin Hochster. 
Michigan Math. J. {\bf{57}} (2008), 249--259. 


\bibitem[EV92]{EV92}
H. Esnault, E. Viehweg. 
\textit{Lectures on vanishing theorems.} 
DMV Seminar, {\bf{20}}. Birkh\"auser Verlag, Basel, (1992). 


\bibitem[Fuj00]{Fuj00}
O. Fujino. 
{\textit{Abundance theorem for semi log canonical threefolds.}} 
Duke Math. J. {\textbf{102}} (2000), no. 3, 513--532. 


\bibitem[Fuj11]{Fuj11}
O. Fujino. 
{\textit{Fundamental theorems for the log minimal model program.}} 
Publ. Res. Inst. Math. Sci. {\textbf{47}} (2011), no. 3, 727--789.


\bibitem[Fuj12a]{Fuj12a}
O. Fujino. 
\textit{A transcendental approach to Koll\'ar's injectivity theorem.}
Osaka J. Math. {\bf{49}} (2012), no. 3, 833--852.


\bibitem[Fuj12b]{Fuj12b}
O. Fujino. 
\textit{A transcendental approach to Koll\'ar's injectivity theorem I\hspace{-.1em}I.}
J. Reine Angew. Math. {\bf{681}} (2013), 149--174. 


\bibitem[Fuj13]{Fuj13}
O. Fujino. 
\textit{Injectivity theorems.}
Preprint, arXiv:1303.2404v1. 



\bibitem[FG14]{FG14}
O. Fujino, Y. Gongyo. 
{\textit{Log pluricanonical representations and the abundance conjecture.}}
Compositio Math.  {\bf{150}} (2014), 593--620. 


\bibitem[GM13]{GM13}
Y. Gongyo, S. Matsumura. 
\textit{Versions of injectivity and extension theorems.}
Preprint, arXiv:1406.6132v2. 


\bibitem[GZ13]{GZ13}
Q. Guan, X. Zhou. 
\textit{Strong openness conjecture for plurisubharmonic functions.}
Preprint, arXiv:1311.3781v1.


\bibitem[Kaw82]{Kaw82}
Y. Kawamata.
\textit{A generalization of Kodaira-Ramanujam's vanishing theorem.}
Math. Ann. {\bf{261}} (1982), no. 1, 43--46.


\bibitem[Kaw85]{Kaw85}
Y. Kawamata. 
\textit {Pluricanonical systems on minimal algebraic varieties.}
Invent. Math. {\bf{79}} (1985), no. 3, 567--588. 


\bibitem[Kaw92]{Kaw92}
Y. Kawamata. 
{\textit{Abundance theorem for minimal threefolds}},   
Invent. Math.  {\textbf{108}}  (1992),  no. 2, 229--246.


\bibitem[KaMM87]{kamama}
Y. Kawamata, K. Matsuda, K. Matsuki.  
\textit{Introduction to the minimal model problem.} 
Algebraic geometry, Sendai, (1985),  283--360, 
Adv. Stud. Pure Math., {\bf{10}}, North-Holland, Amsterdam, (1987). 


\bibitem[KeMMc04]{kemamc}
S. Keel, K. Matsuki, J. $\mathrm{M^{c}}$Kernan.  
{\textit{Log abundance theorem for threefolds.}}  
Duke Math. J.  {\textbf{75}}  (1994),  no. 1, 99--119, 
Corrections to: 
\lq \lq{\textit{Log abundance theorem for threefolds.}}"   
Duke Math. J.  {\textbf{122}}  (2004),  no. 3, 625--630. 

\bibitem[KM]{KM}
J. Koll\'ar, S. Mori. 
{\textit{Birational geometry of algebraic varieties.}}
Cambridge Tracts in Math.,{\bf{134}} (1998).

\bibitem[Kol86]{Kol86}
J. Koll\'ar. 
\textit{Higher direct images of dualizing sheaves. I.}
Ann. of Math. (2) {\bf{123}} (1986), no. 1, 11--42.


\bibitem[Laz]{Laz}
R. Lazarsfeld.  
\textit{Positivity in Algebraic Geometry I-II.}
A Series of Modern Surveys
in Mathematics, {\bf{48, 49}}
Springer Verlag, Berlin, (2004).


\bibitem[Mat13a]{Mat13a}
S. Matsumura. 
\textit{
A Nadel vanishing theorem for metrics with minimal singularities on big line bundles.}
Preprint, arXiv:1306.2497v2.


\bibitem[Mat13b]{Mat13b}
S. Matsumura. 
\textit{An injectivity theorem with multiplier 
ideal sheaves of singular metrics 
with transcendental singularities.}
Preprint, arXiv:1308.2033v2. 


\bibitem[Mat14]{Mat14}
S. Matsumura. 
\textit{A Nadel vanishing theorem via injectivity theorems.}
Math. Ann. {\bf{359}} (2014) no.4, pp. 785-802. 




\bibitem[Nad89]{Nad89}
A. M. Nadel. 
\textit{Multiplier ideal sheaves and existence of 
K\"ahler-Einstein metrics of positive scalar curvature.}
Proc. Nat. Acad. Sci. U.S.A. {\bf{86}} (1989), 
no. 19, 7299--7300. 


\bibitem[Nad90]{Nad90}
A. M. Nadel. 
\textit{Multiplier ideal sheaves and K\"ahler-Einstein metrics of positive scalar curvature.}
Ann. of Math. (2) {\bf{132}} (1990), no. 3, 549--596. 


\bibitem[Nak]{Nak}
N. Nakayama. 
\textit{Zariski decomposition and abundance.} 
MSJ Memoirs, {\bf{14}}. Mathematical Society of Japan, Tokyo, 2004. 


\bibitem[Ohs84]{Ohs84}
T. Ohsawa. 
\textit{Vanishing theorems on complete K\"ahler manifolds.} 
Publ. Res. Inst. Math. Sci. {\bf{20}} (1984), no. 1, 21--38. 


\bibitem[Ohs04]{Ohs04}
T. Ohsawa. 
\textit{On a curvature condition that implies a cohomology injectivity theorem of Koll\'ar-Skoda type.} 
Publ. Res. Inst. Math. Sci. {\bf{41}} (2005), 
no. 3, 565--577.

\bibitem[Tak97]{Tak97}
K. Takegoshi. 
\textit{
On cohomology groups of nef line bundles tensorized with multiplier ideal sheaves on compact K\"ahler manifolds.} 
Osaka J. Math. {\bf 34} (1997), no. 4, 783--802.


\bibitem[Tan71]{Tan71} 
S. G. Tankeev.  
{\textit{On $n$-dimensional canonically polarized varieties 
and varieties of fundamental type.}}
Math. USSR-Izv. {\textbf{5}} (1971), no. 1, 29--43. 



\bibitem[Ver10]{Ver10}
M. Verbitsky.  
{\textit{HyperK$\mathrm{\ddot{a}}$hler SYZ conjecture and semipositive line bundles.}}
Geom. Funct. Anal.  {\bf{19}} (2010), no. 5, 1481--1493.


\bibitem[Vie82]{Vie82}
E. Viehweg. 
\textit{Vanishing theorems.}
J. Reine Angew. Math. {\bf{335}} (1982), 1--8. 



\end{thebibliography}
\end{document}